\documentclass[10pt]{amsart}
\usepackage{amsfonts}
\usepackage{amssymb}
\newtheorem{thm1}{Theorem}[section]
\newtheorem{lem1}[thm1]{Lemma}
\newtheorem{rem1}[thm1]{Remark}
\newtheorem{def1}[thm1]{Definition}

\newtheorem{cor1}[thm1]{Corollary}

\newtheorem{prop1}[thm1]{Proposition}
\newtheorem{ex1}[thm1]{Example}
\newtheorem{not1}[thm1]{Notation}
\begin{document}

\title
{ARITHMETICAL RANK OF BINOMIAL IDEALS}
\author[A. Katsabekis]{Anargyros Katsabekis}
\address {Department of Mathematics, Bilkent University, 06800 Ankara, Turkey} \email{katsampekis@bilkent.edu.tr}

\keywords{Arithmetical rank; binomial ideals; graphs; indispensable monomials.}
\subjclass[2010]{13F20, 14M12, 05C25}

\begin{abstract} In this paper, we investigate the arithmetical rank of a binomial ideal $J$. We provide lower bounds for the binomial arithmetical rank and the $J$-complete arithmetical rank of $J$. Special attention is paid to the case where $J$ is the binomial edge ideal of a graph. We compute the arithmetical rank of such an ideal in various cases.
\end{abstract}
\maketitle


\section{Introduction}

Consider the polynomial ring $K[x_{1},\ldots,x_{m}]$ in the variables $x_{1},\ldots,x_{m}$ over a field $K$. For the sake of simplicity we will denote by ${\bf x}^{\bf u}$ the monomial $x_{1}^{u_1} \cdots x_{m}^{u_m}$ of $K[x_{1},\ldots,x_{m}]$, with ${\bf u}=(u_{1},\ldots,u_{m}) \in \mathbb{N}^{m}$, where $\mathbb{N}$ stands for the set of non-negative integers. A {\em binomial} in the sense of \cite[Chapter 8]{Vil} is a difference of two monomials, i.e. it is of the form ${\bf x}^{\bf u}-{\bf x}^{\bf v}$. A {\em binomial ideal} is an ideal generated by binomials. 

Toric ideals serve as important examples of binomial ideals. Let $\mathcal{A}=\{{\bf a}_{1},\ldots,{\bf a}_{m}\}$ be a subset of $\mathbb{Z}^{n}$. The {\em toric ideal} $I_{\mathcal{A}}$ is the kernel of the $K$-algebra homomorphism $\phi: K[x_{1},\ldots,x_{m}] \rightarrow K[t_{1}^{\pm 1},\ldots,t_{n}^{\pm 1}]$ given by $$\phi(x_{i})= {\bf t}^{{\bf a}_i}=t_{1}^{a_{i,1}} \cdots t_{n}^{a_{i,n}} \ \ \textrm{for all} \ i=1,\ldots,m,$$ where ${\bf a}_{i}=(a_{i,1},\ldots,a_{i,n})$.

We grade $K[x_{1},\ldots,x_{m}]$ by the semigroup $\mathbb{N}\mathcal{A}:=\{l_{1}{\bf a}_{1}+\cdots+l_{m}{\bf a}_{m}|l_{i} \in \mathbb{N}\}$ setting ${\rm deg}_{\mathcal{A}}(x_{i})={\bf a}_{i}$ for $i=1,\ldots,m$. The $\mathcal{A}$-degree of a monomial ${\bf x}^{\bf u}$ is defined by $${\rm deg}_{\mathcal{A}}({\bf x}^{\bf u})=u_{1}{\bf a}_{1}+\cdots+u_{m}{\bf a}_{m} \in \mathbb{N}\mathcal{A}.$$ A polynomial $F \in K[x_{1},\ldots,x_{m}]$ is $\mathcal{A}-$homogeneous if the $\mathcal{A}-$degrees of all the monomials that occur in $F$ are the same. An ideal is $\mathcal{A}$-homogeneous if it is generated by $\mathcal{A}$-homogeneous polynomials. The ideal $I_{\mathcal{A}}$ is generated by all the binomials ${\bf x}^{{\bf u}}-{\bf x}^{{\bf v}}$ such that ${\rm deg}_{\mathcal{A}}({\bf x}^{{\bf u}})={\rm deg}_{\mathcal{A}}({\bf x}^{{\bf v}})$ (see \cite[Lemma 4.1]{St}), thus $I_{\mathcal{A}}$ is $\mathcal{A}$-homogeneous. 

Let $J \subset K[x_{1},\ldots,x_{m}]$ be a binomial ideal. There exist a positive integer $n$ and a vector configuration $\mathcal{A} = \{\mathbf{a}_1, \ldots, \mathbf{a}_m\} \subset \mathbb{Z}^n$ such that $J \subset I_\mathcal{A}$, see for instance \cite[Theorem 1.1]{KaOj}. We say that a polynomial $F=c_{1}M_{1}+\ldots+c_{s}M_{s} \in J$, where $c_{i} \in K$ and $M_{1},\ldots,M_{s}$ are monomials, is $J$-{\em complete} if $M_{i}-M_{l}\in J$ for every $1 \leq i <l \leq s$. Clearly every $J$-complete polynomial $F$ is also $\mathcal{A}$-homogeneous.

Computing the least number of polynomial equations defining an algebraic set is a classical problem in Algebraic Geometry which goes back to Kronecker \cite{Kr}. This problem is equivalent, over an algebraically closed field, with the corresponding problem in Commutative Algebra of the determination of the smallest integer $s$ for which there exist polynomials $F_{1},\ldots,F_{s}$ in $J$ such that $rad(J)=rad(F_{1},\ldots,F_{s})$. The number $s$ is commonly known as the {\em arithmetical rank} of $J$ and will be denoted by ${\rm ara}(J)$. Since $J$ is generated by binomials, it is natural to define the {\em binomial arithmetical rank} of $J$, denoted by ${\rm bar}(J)$, as the smallest integer $s$ for which there exist binomials $B_{1},\ldots,B_{s}$ in $J$ such that $rad(J)=rad(B_{1},\ldots,B_{s})$. Furthermore we can define the $J$-{\em complete arithmetical rank} of $J$, denoted by ${\rm ara}_{c}(J)$, as the smallest integer $s$ for which there exist $J$-complete polynomials $F_{1},\ldots,F_{s}$ in $J$ such that $rad(J)=rad(F_{1},\ldots,F_{s})$. Finally we define the $\mathcal{A}$-{\em homogeneous arithmetical rank} of $J$, denoted by ${\rm ara}_{\mathcal{A}}(J)$, as the smallest integer $s$ for which there exist $\mathcal{A}$-homogeneous polynomials $F_{1},\ldots,F_{s}$ in $J$ such that $rad(J)=rad(F_{1},\ldots,F_{s})$. From the definitions and \cite[Corollary 3.3.3]{BS98} we deduce the following inequalities: $${\rm cd}(J) \leq {\rm ara}(J) \leq {\rm ara}_{\mathcal{A}}(J)\leq {\rm ara}_{c}(J) \leq {\rm bar}(J)$$ where ${\rm cd}(J)$ is the cohomological dimension of $J$.

In section 2 we introduce the simplicial complex $\Delta_{J}$ and use combinatorial invariants of the aforementioned complex to provide lower bounds for the binomial arithmetical rank and the $J$-complete arithmetical rank of $J$. In particular we prove that ${\rm bar}(J) \geq \delta(\Delta_{J})_{\{0,1\}}$ and ${\rm ara}_{c}(J) \geq \delta(\Delta_{J})_{\Omega}$, see Theorem \ref{Basic}.

In section 3 we study the arithmetical rank of the binomial edge ideal $J_{G}$ of a graph $G$. This class of ideals generalize naturally the determinantal ideal generated by the 2-minors of the matrix $$\begin{pmatrix}
x_{1} & x_{2} & \ldots & x_{n}\\
x_{n+1} & x_{n+2} & \ldots & x_{2n}\end{pmatrix}.$$We prove (see Theorem \ref{Basbar}) that, for a binomial edge ideal $J_{G}$, both the binomial arithmetical rank and the $J_{G}$-complete arithmetical rank coincide with the number of edges of $G$. If $G$ is the complete graph on the vertex set $\{1,\ldots,n\}$, then, from \cite[Theorem 2]{BruS}, the arithmetical rank of $J_{G}$ equals $2n-3$. It is still an open problem to compute ${\rm ara}(J_{G})$ when $G$ is not the complete graph. We show that ${\rm ara}(J_{G}) \geq n+l-2$, where $n$ is the number of vertices of $G$ and $l$ is the vertex connectivity of $G$. Furthermore we prove that in several cases ${\rm ara}(J_{G})={\rm cd}(J_{G})=n+l-2$, see Theorem \ref{Cycle}, Theorem \ref{Exactlyone} and Theorem \ref{many odd}.



\section{Lower bounds}

First we will use the notion of indispensability to introduce the simplicial complex $\Delta_{J}$. Let $J \subset K[x_{1},\ldots,x_{m}]$ be a binomial ideal containing no binomials of the form ${\bf x}^{\bf u}-1$, where ${\bf u} \neq {\bf 0}$. A {\em binomial} $B=M-N \in J$ is called {\em indispensable} of $J$ if every system of binomial generators of $J$ contains $B$ or $-B,$ while a {\em monomial} $M$ is called {\em indispensable} of $J$ if every system of binomial generators of $J$ contains a binomial $B$ such that $M$ is a monomial of $B$. Let $\mathcal{M}_{J}$ be the ideal generated by all monomials $M$ for which there exists a nonzero $M-N \in J$. By \cite[Proposition 1.5]{KaOj} the set $G(\mathcal{M}_{J})$ of indispensable monomials of $J$ is the unique minimal generating set of $\mathcal{M}_{J}$.

The support of a monomial ${\bf x}^{\bf u}$ of $K[x_{1},\ldots,x_{m}]$ is ${\rm supp}({\bf x}^{\bf u}):=\{i |x_{i} \ \textrm{divides} \ {\bf x}^{\bf u}\}$. Let $\mathcal{T}$ be the set of all $E \subset \{1,\ldots,m\}$ for which there exists an indispensable monomial $M$ of $J$ such that $E={\rm supp}(M)$. Let $\mathcal{T}_{\rm min}$ denote the set of minimal elements of $\mathcal{T}$.

\begin{def1} {\rm  We associate to $J$ a simplicial complex $\Delta_{J}$ with vertices the elements of $\mathcal{T}_{\rm min}$. Let $T=\{E_{1},\ldots,E_{k}\}$ be a subset of $\mathcal{T}_{\rm min}$, then $T \in \Delta_{J}$ if there exist $M_{i}$, $1 \leq i \leq k$, such that ${\rm supp}(M_{i})=E_{i}$ and $M_{i}-M_{l} \in J$ for every $1 \leq i<l \leq k$.}
\end{def1}

Next we will study the connection between the radical of $J$ and $\Delta_{J}$. The {\em induced subcomplex} $\Delta'$ of $\Delta_{J}$ by certain vertices $\mathcal{V} \subset \mathcal{T}_{\rm min}$ is the subcomplex of $\Delta_{J}$ with vertices $\mathcal{V}$ and $T \subset \mathcal{V}$ is a simplex of the subcomplex $\Delta'$ if $T$ is a simplex of $\Delta_{J}$. A subcomplex $H$ of $\Delta_{J}$ is called a {\em spanning subcomplex} if both have exactly the same set of vertices.

Let $F$ be a polynomial in $K[x_{1},\ldots,x_{m}]$. We associate to $F$ the induced subcomplex $\Delta_{J}(F)$ of $\Delta_{J}$ consisting of those vertices $E_{i} \in \mathcal{T}_{\rm min}$ with the property: there exists a monomial $M_{i}$ in $F$ such that $E_{i}={\rm supp}(M_{i})$. The next theorem provides a necessary condition under which a set of polynomials in the binomial ideal $J$ generates the radical of $J$ up to radical.

\begin{prop1} \label{Spanning} Let $K$ be any field. If $rad(J)=rad(F_{1},\ldots,F_{s})$ for some polynomials $F_{1},\ldots,F_{s}$ in $J$, then $\cup_{i=1}^{s}\Delta_{J}(F_{i})$ is a spanning subcomplex of $\Delta_{J}$.
\end{prop1}
\noindent \noindent \textbf{Proof.} Let $E={\rm supp}({\bf x}^{{\bf u}}) \in \mathcal{T}_{\rm min}$, where $B={\bf x}^{{\bf u}}-{\bf x}^{{\bf v}} \in J$ and ${\bf x}^{{\bf u}}$ is an indispensable monomial of $J$. We will show that there exists a monomial $M$ in some $F_{l}$, $1 \leq l \leq s$, such that $E={\rm supp}(M)$. Since $rad(J)=rad(F_{1},\ldots,F_{s})$, there is a power $B^{r}$, $r \geq 1$, which belongs to the ideal generated by $F_{1},\ldots,F_{s}$. Thus there is a monomial $M$ in some $F_{l}$ dividing the monomial $({\bf x}^{{\bf u}})^{r}$ and therefore ${\rm supp}(M) \subseteq {\rm supp}({\bf x}^{{\bf u}})$. But $F_{l} \in J$ and $J$ is generated by binomials, so there exists ${\bf x}^{{\bf z}}-{\bf x}^{{\bf w}} \in J$ such that ${\bf x}^{{\bf z}}$ divides $M$. Since ${\bf x}^{\bf z} \in \mathcal{M}_{J}$ and $G(\mathcal{M}_{J})$ generates $\mathcal{M}_{J}$, there is an indispensable monomial $N$ dividing ${\bf x}^{{\bf z}}$, thus $${\rm supp}(N) \subseteq {\rm supp}({\bf x}^{\bf z}) \subseteq {\rm supp}(M) \subseteq E.$$  Since $E \in \mathcal{T}_{\rm min}$, we deduce that $E={\rm supp}(N)$, and therefore $E={\rm supp}(M)$. \hfill $\square$

\begin{rem1} \label{RemA} {\rm (1) If $F$ is a $J$-complete polynomial of $J$, then $\Delta_{J}(F)$ is a simplex. To see that $\Delta_{J}(F)$ is a simplex suppose that $\Delta_{J}(F) \neq \emptyset$ and let $T=\{E_{1},\ldots,E_{k}\}$ be the set of vertices of $\Delta_{J}(F)$. For every $1 \leq i \leq k$ there exists a monomial $M_{i}$, $1 \leq i \leq k$, in $F$ such that $E_{i}={\rm supp}(M_{i})$. Since $F$ is $J$-complete, we have that $M_{i}-M_{l} \in J$ for every $1 \leq i<l \leq k$. Thus $\Delta_{J}(F)$ is a simplex.\\ (2) If $B$ is a binomial of $J$, then $\Delta_{J}(B)$ is either a vertex, an edge or the empty set.}
\end{rem1}

\begin{rem1} \label{Basic-Bar} {\rm If the equality $rad(J)=rad(F_{1},\ldots,F_{s})$ holds for some $J$-complete polynomials $F_{1},\ldots,F_{s}$ in $J$, then $\cup_{i=1}^{s} \Delta_{J}(F_{i})$ is a spanning subcomplex of $\Delta_{J}$ and each $\Delta_{J}(F_{i})$ is a simplex.}
\end{rem1}

For a simplicial complex $\Delta$ we denote by $r_{\Delta}$ the smallest number $s$ of simplices $T_{i}$ of $\Delta$, such that the subcomplex $\cup_{i=1}^{s}T_{i}$ is spanning and by $b_{\Delta}$ the smallest number $s$ of simplices $T_{i}$ of $\Delta$, such that the subcomplex $\cup_{i=1}^{s}T_{i}$ is spanning and each $T_{i}$ is either an edge, a vertex or the empty set.

\begin{thm1} \label{Basic1} Let $K$ be any field, then $b_{\Delta_{J}} \leq {\rm bar}(J)$ and $r_{\Delta_{J}} \leq {\rm ara}_{c}(J)$.
\end{thm1}

It turns out that both $b_{\Delta_{J}}$ and $r_{\Delta_{J}}$ have a combinatorial interpretation in terms of matchings in $\Delta_{J}$.

Let $\Delta$ be a simplicial complex on the vertex set $\mathcal{T}_{\rm min}$ and $Q$ be a subset of $\Omega:=\{0,1,\ldots,{\rm dim}(\Delta)\}$. A set $\mathcal{N}=\{T_{1},\ldots,T_{s}\}$ of simplices of $\Delta$ is called a $Q$-{\em matching} in $\Delta$ if $T_{k} \cap T_{l}=\emptyset$ for every $1 \leq k,l \leq s$ and ${\rm dim}(T_{k}) \in Q$ for every $1 \leq k \leq s$; see also Definition 2.1 in \cite{KT}. Let ${\rm supp}(\mathcal{N})=\cup_{i=1}^{s}T_{i}$, which is a subset of the vertices $\mathcal{T}_{\rm min}$. We denote by ${\rm card}(\mathcal{N})$ the cardinality $s$ of the set $\mathcal{N}$. A $Q$-matching $\mathcal{N}$ in $\Delta$ is called a {\em maximal} $Q$-{\em matching} if ${\rm supp}(\mathcal{N})$ has the maximum possible cardinality among all $Q$-matchings. By $\delta(\Delta)_{Q}$, we denote the minimum of the set $$\{{\rm card}(\mathcal{N})| \mathcal{N} \ \textrm{is a maximal} \ Q-\textrm{matching in} \ \Delta\}.$$

\begin{thm1} \label{Basic} Let $K$ be any field, then ${\rm bar}(J) \geq \delta(\Delta_{J})_{\{0,1\}}$ and ${\rm ara}_{c}(J) \geq \delta(\Delta_{J})_{\Omega}$.
\end{thm1}
\noindent \noindent \textbf{Proof.} By Proposition 3.3 in \cite{KT}, $b_{\Delta_{J}}=\delta(\Delta_{J})_{\{0,1\}}$ and $r_{\Delta_{J}}=\delta(\Delta_{J})_{\Omega}$. Now the result follows from Theorem \ref{Basic1}. \hfill $\square$

\begin{prop1} \label{Indispenbasic} Let $J$ be a binomial ideal. Suppose that there exists a minimal generating set $\mathcal{S}$ of $J$ such that every element of $\mathcal{S}$ is a difference of two squarefree monomials. Assume that $J$ is generated by the indispensable binomials, namely $\mathcal{S}$ consists precisely of the indispensable binomials (up to sign). Then ${\rm bar}(J)={\rm card}(\mathcal{S})$.
\end{prop1}
\noindent \textbf{Proof.}  Let ${\rm card}(\mathcal{S})=t$. Since $\mathcal{S}$ is a generating set of $J$, we have that ${\rm bar}(J) \leq t$. It is enough to prove that $t \leq {\rm bar}(J)$. Let $|\mathcal{T}_{\rm min}|=g$. By \cite[Corollary 3.6]{CTV} it holds that ${\rm card}(G(\mathcal{M}_{J}))=2t$, so $g=2t$. For every maximal $\{0,1\}$-matching $\mathcal{M}$ in $\Delta_{J}$ we have that ${\rm supp}(\mathcal{M})=\mathcal{T}_{min}$, so $\delta(\Delta_{J})_{\{0,1\}} \geq \left \lfloor{\frac{g}{2}} \right \rfloor$ and therefore $\delta(\Delta_{J})_{\{0,1\}}  \geq t$. Thus, from Theorem \ref{Basic}, ${\rm bar}(J) \geq t$. \hfill $\square$

\begin{ex1} {\rm Let $J$ be the binomial ideal generated by $f_{1}=x_{1}x_{6}-x_{2}x_{5}$, $f_{2}=x_{2}x_{7}-x_{3}x_{6}$, $f_{3}=x_{1}x_{8}-x_{4}x_{5}$, $f_{4}=x_{3}x_{8}-x_{4}x_{7}$ and $f_{5}=x_{1}x_{7}-x_{3}x_{5}$. Actually $J$ is the binomial edge ideal of the graph $G$ with edges $\{1,2\}$, $\{2,3\}$, $\{1,4\}$, $\{3,4\}$ and $\{1,3\}$, see section 3 for the definition of such an ideal. Note that $J$ is $\mathcal{A}$-homogeneous where $\mathcal{A}=\{{\bf a}_{1},\ldots,{\bf a}_{8}\}$ is the set of columns of the matrix $$D=\begin{pmatrix}
1 & 0 & 0 & 0 & 1 & 0 & 0 & 0 
            \\
           0 & 1 & 0 & 0 & 0 & 1 & 0 & 0 
            \\
            0 & 0 & 1 & 0 & 0 & 0 & 1 & 0 
            \\
            0 & 0 & 0 & 1 & 0 & 0 & 0 & 1

            \end{pmatrix}.$$

By \cite[Theorem 3.3]{CTV} every binomial $f_{i}$ is indispensable of $J$. Thus $$\mathcal{T}_{\rm min}=\{E_{1}=\{1,6\}, E_{2}=\{2,5\}, E_{3}=\{2,7\}, E_{4}=\{3,6\}, E_{5}=\{1,8\},$$ $$E_{6}=\{4,5\}, E_{7}=\{3,8\}, E_{8}=\{4,7\}, E_{9}=\{1,7\}, E_{10}=\{3,5\}\}.$$ By Proposition \ref{Indispenbasic} the binomial arithmetical rank of $J$ equals 5. The simplicial complex $\Delta_{J}$ has 5 connected components and all of them are 1-simplices, namely $\Delta_{1}=\{E_{1},E_{2}\}$, $\Delta_{2}=\{E_{3},E_{4}\}$, $\Delta_{3}=\{E_{5},E_{6}\}$, $\Delta_{4}=\{E_{7},E_{8}\}$ and $\Delta_{5}=\{E_{9},E_{10}\}$. Consequently $$\delta(\Delta_{J})_{\Omega}=\sum_{i=1}^{5}\delta(\Delta_{i})_{\Omega}=1+1+1+1+1=5$$ and therefore $5 \leq {\rm ara}_{c}(J)$. Since ${\rm ara}_{c}(J) \leq {\rm bar}(J)$, we get that ${\rm ara}_{c}(J)=5$. We will show that ${\rm ara}_{\mathcal{A}}(J)=5$. Suppose that ${\rm ara}_{\mathcal{A}}(J)=s<5$ and let $F_{1},\ldots,F_{s}$ be $\mathcal{A}$-homogeneous polynomials in $J$ such that $rad(J)=rad(F_{1},\ldots,F_{s})$. For every vertex $E_{i} \in \mathcal{T}_{\rm min}$ there exists, from Proposition \ref{Spanning}, a monomial $M_{i}$ in $F_{k}$ such that $E_{i}={\rm supp}(M_{i})$. But $s<5$, so there exist $E_{i} \in \mathcal{T}_{\rm min}$ and $E_{j} \in \mathcal{T}_{\rm min}$ such that \begin{enumerate} \item $\{E_{i}, E_{j}\}$ is not a $1$-simplex of $\Delta_{J}$,\item $E_{i}={\rm supp}(M_{i})$, $E_{j}={\rm supp}(M_{j})$ and \item $M_{i}$ and $M_{j}$ are monomials of some $F_{k}$.
\end{enumerate}
Since $F_{k}$ is $\mathcal{A}$-homogeneous, it holds that ${\rm deg}_{\mathcal{A}}(M_{i})={\rm deg}_{\mathcal{A}}(M_{j})$. Considering all possible combinations of $E_{i}$ and $E_{j}$ we finally arrive at a contradiction. Thus ${\rm ara}_{\mathcal{A}}(J)=5$. Note that $J$ is $\mathcal{B}$-homogeneous where $\mathcal{B}$ is the set of columns of the matrix $$N=\begin{pmatrix}
1 & 1 & 1 & 1 & 0 & 0 & 0 & 0 
            \\
           1 & 0 & 0 & 0 & 1 & 0 & 0 & 0 
            \\
            0 & 1 & 0 & 0 & 0 & 1 & 0 & 0 
            \\
            0 & 0 & 1 & 0 & 0 & 0 & 1 & 0\\
0 & 0 & 0 & 1 & 0 & 0 & 0 & 1
            \end{pmatrix}.$$ Since every row of $D$ is a row of $N$, we deduce that every $\mathcal{B}$-homogeneous polynomial in $J$ is also $\mathcal{A}$-homogeneous. So ${\rm ara}_{\mathcal{B}}(J)$ is an upper bound for ${\rm ara}_{\mathcal{A}}(J)$, therefore ${\rm ara}_{\mathcal{B}}(J)=5$. We have that $rad(J)=rad(f_{1},f_{2}+f_{3},f_{4},f_{5})$, since the second power of both binomials $f_{2}$ and $f_{3}$ belongs to the ideal generated by the polynomials $f_{1},f_{2}+f_{3},f_{4},f_{5}$. Remark that the polynomials $f_{1}$, $f_{2}+f_{3}$, $f_{4}$ and $f_{5}$ are $\mathcal{C}$-homogeneous, where $\mathcal{C}$ is the set of columns of the matrix $$\begin{pmatrix}
1 & 2 & 3 & 4 & 1 & 2 & 3 & 4 
            \\
            5 & 6 & 7 & 8 & 5 & 6 & 7 & 8 

            \end{pmatrix}.$$ Thus ${\rm ara}_{\mathcal{C}}(J) \leq 4$, so ${\rm ara}(J) \leq 4$. A primary decomposition of $J$ is $$J=(f_{1},f_{2},f_{3},f_{4},f_{5},x_{2}x_{8}-x_{4}x_{6})\cap (x_{1},x_{3},x_{5},x_{7}).$$ Hence, by \cite[Proposition 19.2.7]{BS98}, it follows that ${\rm ara}(J) \geq 4$. Thus $${\rm ara}(J)={\rm ara}_{\mathcal{C}}(J)=4<5={\rm ara}_{\mathcal{A}}(J)={\rm ara}_{\mathcal{B}}(J)={\rm ara}_{c}(J)={\rm bar}(J).$$}
\end{ex1}



\section{Binomial edge ideals of graphs}

In this section we consider a special class of binomial ideals, namely binomial edge ideals of graphs. This ideal was introduced in \cite{HTH} and independently at the same time in \cite{Ohtani}.

Let $G$ be an undirected connected simple graph on the vertex set $[n]:=\{1,\ldots,n\}$ and with edge set $E(G)$. Consider the polynomial ring $$R:=K[x_{1},\ldots,x_{n},x_{n+1},\ldots,x_{2n}]$$ in $2n$ variables, $x_1, \ldots, x_n,$ $x_{n+1},\ldots,x_{2n}$, over $K$.

\begin{def1}
{\rm The binomial edge ideal $J_G \subset R$ associated to the graph $G$ is the ideal generated by the binomials $f_{ij} = x_ix_{n+j}- x_{j}x_{n+i},$ with $i < j,$ such that $\{i, j\}$ is an edge of $G.$}
\end{def1}

\begin{rem1}  {\rm From \cite[Corollary 1.13]{KaOj} every binomial $f_{ij}$, where $\{i,j\}$ is an edge of $G$, is indispensable of $J_{G}$. Thus $$\mathcal{T}_{min}=\{E_{ij}^{1}=\{i,n+j\},E_{ij}^{2}=\{j,n+i\}|\{i,j\} \in E(G)\}.$$}
\end{rem1}

We recall some fundamental material from \cite{HTH}. Let $G$ be a connected graph on $[n]$ and let $S \subset [n]$. By $G \setminus S$, we denote the graph that results from deleting all vertices in $S$ and their incident edges from $G$. Let $c(S)$ be the number of connected components of $G \setminus S$ and let $G_{1},\ldots,G_{c(S)}$ denote the connected components of $G \setminus S$. Also let $\overset{\sim}{G_{i}}$ denote the complete graph on the vertices of $G_{i}$. We set $$P_{S}(G)=(\cup_{i \in S} \{x_{i},x_{n+i}\}, J_{\overset{\sim}{G_{1}}},\ldots,J_{\overset{\sim}{G}_{c(S)}})R.$$ Then $P_{S}(G)$ is a prime ideal for every $S \subset [n]$. The ring $R/P_{\emptyset}(G)$ has Krull dimension $n+1$. For $S \neq \emptyset$ the ring $R/P_{S}(G)$ has Krull dimension $n-{\rm card}(S)+c(S)$. The ideal $P_{S}(G)$ is a minimal prime of $J_{G}$ if and only if $S=\emptyset$ or $S \neq \emptyset$ and for each $i \in S$ one has $c(S \setminus \{i\})<c(S)$. Moreover $J_{G}$ is a radical ideal and it admits the minimal primary decomposition $J_{G}= \cap_{S \in \mathcal{M}(G)}P_{S}(G)$, where $\mathcal{M}(G)=\{S \subset [n]: P_{S}(G) \ \textrm{is a minimal prime of} \ J_{G}\}$.

\begin{thm1} \label{Basbar} Let $G$ be a connected graph on the vertex set $[n]$ with $m$ edges. Then ${\rm bar}(J_{G})={\rm ara}_{c}(J_{G})=m$.

\end{thm1}

\noindent \textbf{Proof.} Every binomial $f_{ij}$, where $\{i,j\}$ is an edge of $G$, is indispensable of $J_{G}$, thus, from Proposition \ref{Indispenbasic}, ${\rm bar}(J_{G})=m$. Note that, for every edge $\{i,j\}$ of $G$, $\{E_{ij}^{1},E_{ij}^{2}\}$ is a 1-simplex of $\Delta_{J_{G}}$. Furthermore $\Delta_{J_{G}}$ has exactly $m$ connected components and all of them are 1-simplices. Thus $\delta(\Delta_{J_{G}})_{\Omega}=m$ and therefore, from Theorem \ref{Basic}, ${\rm ara}_{c}(J_{G}) \geq m$. Consequently ${\rm ara}_{c}(J_{G})=m$. \hfill $\square$

\begin{thm1} \label{Propara} Let $G$ be a connected graph on the vertex set $[n]$ with $m$ edges. Consider the canonical basis $\{{\bf e}_{1},\ldots,{\bf e}_{n}\}$ of $\mathbb{Z}^{n}$ and the canonical basis $\{{\bf w}_{1},\ldots,{\bf w}_{n+1}\}$ of $\mathbb{Z}^{n+1}$. Let $\mathcal{A}=\{{\bf a}_{1}, \ldots,{\bf a}_{2n}\} \subset \mathbb{N}^{n}$ be the set of vectors where ${\bf a}_{i}={\bf e}_{i}$, $1 \leq i \leq n$, and ${\bf a}_{n+i}={\bf e}_{i}$ for $1 \leq i \leq n$. Let $\mathcal{B}=\{{\bf b}_{1}, \ldots,{\bf b}_{2n}\} \subset \mathbb{N}^{n+1}$ be the set of vectors where ${\bf b}_{i}={\bf w}_{1}+{\bf w}_{i+1}$, $1 \leq i \leq n$, and ${\bf b}_{n+i}={\bf w}_{i+1}$ for $1 \leq i \leq n$. Then ${\rm ara}_{\mathcal{A}}(J_{G})={\rm ara}_{\mathcal{B}}(J_{G})=m$.

\end{thm1}

\noindent \textbf{Proof.} Suppose that ${\rm ara}_{\mathcal{A}}(J_{G})=t<m$ and let $F_{1},\ldots,F_{t}$ be $\mathcal{A}$-homogeneous polynomials in $J_{G}$ such that $J_{G}=rad(F_{1},\ldots,F_{t})$. For every edge $\{i,j\}$ of $G$ with $i<j$ there exist, from Proposition \ref{Spanning}, monomials $M_{ij}^{k}$ and $N_{ij}^{l}$ in $F_{k}$ and $F_{l}$, respectively, such that $E_{ij}^{1}={\rm supp}(M_{ij}^{k})$ and $E_{ij}^{2}={\rm supp}(N_{ij}^{l})$. But $t<m$, so there exists $E_{rs}^{1} \in \mathcal{T}_{\rm min}$, where $\{r,s\}$ is an edge of $G$ with $r<s$, such that \begin{enumerate} \item $\{E_{ij}^{1}, E_{rs}^{1}\}$ is not a $1$-simplex of $\Delta_{J_{G}}$,\item $E_{ij}^{1}={\rm supp}(M_{ij}^{k})$, $E_{rs}^{1}={\rm supp}(M_{rs}^{k})$ and \item $M_{ij}^{k}$ and $M_{rs}^{k}$ are monomials of some $F_{k}$. 
\end{enumerate}
Let $M_{ij}^{k}=x_{i}^{g_i}x_{n+j}^{g_j}$ and $M_{rs}^{k}=x_{r}^{g_r}x_{n+s}^{g_s}$. Since $F_{k}$ is $\mathcal{A}$-homogeneous, we deduce that ${\rm deg}_{\mathcal{A}}(M_{ij}^{k})={\rm deg}_{\mathcal{A}}(M_{rs}^{k})$, and therefore $g_{i}{\bf e}_{i}+g_{j}{\bf e}_{j}=g_{r}{\bf e}_{r}+g_{s}{\bf e}_{s}$. Consequently $i=r$, $j=s$ and also $M_{ij}^{k}=M_{rs}^{k}$, a contradiction. Let $D$ and $Q$ be the matrices with columns $\mathcal{A}$ and $\mathcal{B}$, respectively. Since every row of $D$ is a row of $Q$, we deduce that every $\mathcal{B}$-homogeneous polynomial in $J_{G}$ is also $\mathcal{A}$-homogeneous. Thus ${\rm ara}_{\mathcal{B}}(J_{G})$ is an upper bound for ${\rm ara}_{\mathcal{A}}(J_{G})$, so $m \leq {\rm ara}_{\mathcal{B}}(J_{G})$ and therefore ${\rm ara}_{\mathcal{B}}(J_{G})=m$. \hfill $\square$

The graph $G$ is called $l$-{\em vertex-connected} if $l<n$ and $G \setminus S$ is connected for every subset $S$ of $[n]$ with ${\rm card}(S)<l$. The {\em vertex connectivity} of $G$ is defined as the maximum integer $l$ such that $G$ is $l$-vertex-connected.

In \cite{BBE} the authors study the relationship between algebraic properties of a binomial edge ideal $J_{G}$, such as the dimension and the depth of $R/J_{G}$, and the vertex connectivity of the graph. It turns out that this notion is also useful for the computation of the arithmetical rank of a binomial edge ideal.

\begin{thm1} \label{Vertexconnectivity} Let $K$ be a field of any characteristic and $G$ be a connected graph on the vertex set $[n]$. Suppose that the vertex connectivity of $G$ is $l$. Then ${\rm ara}(J_{G}) \geq n+l-2$.

\end{thm1}

\noindent \textbf{Proof.} If $G$ is the complete graph on the vertex set $[n]$, its vertex connectivity is $n-1$, then ${\rm ara}(J_{G})=2n-3=n+l-2$ by \cite[Theorem 2]{BruS}. Assume now that $G$ is not the complete graph.  Let $P_{\emptyset}(G)$, $W_{1},\ldots,W_{t}$ be the minimal primes of $J_{G}$. It holds that $J_{G}=P_{\emptyset}(G) \cap L$ where $L=\cap_{i=1}^{t}W_{i}$. First we will prove that ${\rm dim} \left( R/(P_{\emptyset}(G)+L) \right) \leq n-l+1$. For every prime ideal $Q$ such that $P_{\emptyset}(G)+L \subseteq Q$, we have that $L \subseteq Q$, so there is $1 \leq i \leq t$ such that $W_{i} \subseteq Q$. Thus $P_{\emptyset}(G)+W_{i} \subseteq Q$ and therefore ${\rm dim} \left( R/(P_{\emptyset}(G)+L) \right) \leq {\rm dim} \left( R/(P_{\emptyset}(G)+W_{i}) \right)$. It is enough to show that ${\rm dim} \left( R/(P_{\emptyset}(G)+W_{i}) \right) \leq n-l+1$. Let $W_{i}=P_{S}(G)$ for $\emptyset \neq S \subset [n]$. We have that $P_{\emptyset}(G)+P_{S}(G)$ is generated by $$\{x_{i}x_{n+j}-x_{j}x_{n+i}: i, j \in [n]\setminus S \} \cup \{x_{i}, x_{n+i}: i \in S\}.$$ Then ${\rm dim} \left( R/(P_{\emptyset}(G)+P_{S}(G)) \right)=n-{\rm card}(S)+1$. If $l=1$, then ${\rm card}(S) \geq 1$ since $S \neq \emptyset$, and therefore ${\rm dim} \left( R/(P_{\emptyset}(G)+W_{i}) \right) \leq n$. Suppose that $l \geq 2$ and also that ${\rm card}(S)<l$. Since $P_{S}(G)$ is a minimal prime, for every $i \in S$ we have that $c(S \setminus \{i\})<c(S)$. But $G$ is $l$-vertex-connected, namely $G \setminus S$ is connected, so $P_{\emptyset}(G) \subset P_{S}(G)$ a contradiction to the fact that $P_{S}(G)$ is a minimal prime. Thus ${\rm dim} \left( R/(P_{\emptyset}(G)+W_{i}) \right) \leq n-l+1$ and therefore ${\rm dim} \left( R/(P_{\emptyset}(G)+L) \right) \leq n-l+1$. Next we will show that ${\rm min}\{{\rm dim} \left( R/P_{\emptyset}(G) \right), {\rm dim} \left( R/L \right)\}>{\rm dim} \left( R/(P_{\emptyset}(G)+L) \right)$. Recall that ${\rm dim} \left( R/P_{\emptyset}(G) \right)=n+1$, so ${\rm dim} \left( R/(P_{\emptyset}(G)+L) \right)<{\rm dim} \left( R/P_{\emptyset}(G) \right)$. Since $L \subset P_{\emptyset}(G)+L$, we deduce that ${\rm dim} \left( R/(P_{\emptyset}(G)+L) \right) \leq {\rm dim} \left( R/L \right)$. Suppose that ${\rm dim} \left( R/(P_{\emptyset}(G)+L) \right)={\rm dim} \left( R/L \right)$, say equal to $s$, and let $Q_{1} \subsetneqq Q_{2} \subsetneqq \cdots \subsetneqq Q_{s}$ be a chain of prime ideals containing $P_{\emptyset}(G)+L$. Then there is $1 \leq j \leq t$ such that $Q_{1}=W_{j}$. So $P_{\emptyset}(G) \subset W_{j}$, a contradiction. By \cite[Proposition 19.2.7]{BS98} it holds that $${\rm cd}(J_{G}) \geq {\rm dim}(R)-{\rm dim} \left( R/(P_{\emptyset}(G)+L) \right)-1=2n-{\rm dim} \left( R/(P_{\emptyset}(G)+L) \right)-1 \geq$$ $$2n-(n-l+1)-1=n+l-2.$$ Consequently ${\rm ara}(J_{G}) \geq n+l-2$. $\hfill \square$

\begin{ex1} {\rm Let $G$ be the graph on the vertex set $[5]$ with edges $\{1,2\}$, $\{2,3\}$, $\{1,3\}$, $\{2,4\}$, $\{4,5\}$ and $\{3,5\}$. Here the vertex connectivity is $l=2$. By Theorem \ref{Vertexconnectivity}, ${\rm ara}(J_{G}) \geq 5$. The ideal $J_{G}$ is generated up to radical by the polynomials $f_{12},f_{23},f_{13}+f_{24},f_{35}$ and $f_{45}$, since both $f_{13}^{2}$ and $f_{24}^{2}$ belong to the ideal generated by $f_{12},f_{23},f_{13}+f_{24},f_{35}$ and $f_{45}$. Thus ${\rm ara}(J_{G})=5<6={\rm bar}(J_{G})$.}

\end{ex1}

\begin{thm1} \label{Cycle} If $G$ is a cycle of length $n \geq 3$, then ${\rm ara}(J_{G})={\rm bar}(J_{G})=n$.

\end{thm1}

\noindent \textbf{Proof.} The vertex connectivity of $G$ is $2$, so, from Theorem \ref{Vertexconnectivity}, the inequality $n \leq {\rm ara}(J_{G})$ holds. Since $G$ has $n$ edges, we have  that ${\rm ara}(J_{G}) \leq {\rm bar}(J_{G})=n$ and therefore ${\rm ara}(J_{G})=n$. \hfill $\square$

\begin{prop1} \label{Basicara} Let $G$ be a connected graph on $[n]$, with $m$ edges and $n \geq 4$. If $G$ contains an odd cycle of length 3, then ${\rm ara}(J_{G}) \leq m-1$.

\end{prop1}

\noindent \textbf{Proof.} Let $C$ be an odd cycle of $G$ of length 3, with edge set $\{\{1,2\},\{2,3\},\{1,3\}\}$. Since $G$ is connected, without loss of generality there is a vertex $4 \leq i \leq n$ such that $\{1,i\}$ is an edge of $G$. We will show that $(x_{1}x_{n+i}-x_{i}x_{n+1})^{2}$ belongs to the ideal $L$ generated by the polynomials $f_{12},f_{13},f_{1i}+f_{23}$. We have that $$x_{1}^{2}x_{n+i}^{2} \equiv x_{1}x_{n+i}x_{i}x_{n+1}-x_{1}x_{2}x_{n+i}x_{n+3}+x_{1}x_{3} x_{n+i}x_{n+2} \equiv$$ $$x_{1}x_{i}x_{n+i}x_{n+1}-x_{2}x_{n+i}x_{3}x_{n+1}+x_{2}x_{3}x_{n+1}x_{n+i} \equiv x_{1}x_{i}x_{n+i}x_{n+1} \ {\rm mod} \ L.$$ Similarly we have that $x_{i}^{2}x_{n+1}^{2} \equiv x_{1}x_{i}x_{n+i}x_{n+1} \ {\rm mod} \ L$. Thus $x_{1}^{2}x_{n+i}^{2}+x_{i}^{2}x_{n+1}^{2} \equiv 2 x_{1}x_{i}x_{n+i}x_{n+1} \ {\rm mod} \ L$, so $(x_{1}x_{n+i}-x_{i}x_{n+1})^{2}$ belongs to $L$. Next we prove that $(x_{2}x_{n+3}-x_{3}x_{n+2})^2$ belongs to $L$. We have that $$x_{2}^{2}x_{n+3}^{2} \equiv x_{2}x_{n+3}x_{3}x_{n+2}-x_{2}x_{n+3}x_{1}x_{n+i}+x_{2}x_{n+3}x_{i}x_{n+1} \equiv$$ $$x_{2}x_{n+3}x_{3}x_{n+2}-x_{2}x_{n+i}x_{3}x_{n+1}+x_{n+3}x_{i}x_{1}x_{n+2} \equiv$$ $$x_{2}x_{n+3}x_{3}x_{n+2}-x_{1}x_{n+2}x_{n+i}x_{3}+x_{i}x_{n+2}x_{3}x_{n+1} \ {\rm mod} \ L.$$ Furthermore $$x_{3}^{2}x_{n+2}^{2} \equiv x_{2}x_{n+3}x_{3}x_{n+2}-x_{3}x_{n+2}x_{i}x_{n+1}+x_{3}x_{n+2}x_{1}x_{n+i} \ {\rm mod} \ L.$$ Thus $x_{2}^{2}x_{n+3}^{2}+x_{3}^{2}x_{n+2}^2 \equiv 2x_{2}x_{n+3}x_{3}x_{n+2}  \ {\rm mod} \ L$, so $(x_{2}x_{n+3}-x_{3}x_{n+2})^2 \in L$. Let $H$ be the subgraph of $G$ consisting of the cycle $C$ and the edge $\{1,i\}$. Then $J_{G}$ is generated up to radical by the following set of $m-1$ binomials: $$\{f_{kl}| \{k,l\} \in E(G) \setminus E(H)\} \cup \{f_{12},f_{13},f_{1i}+f_{23}\}.$$ Therefore ${\rm ara}(J_{G}) \leq m-1$. \hfill $\square$\\

Let $G_{1}=(V(G_{1}),E(G_{1}))$, $G_{2}=(V(G_{2}),E(G_{2}))$ be graphs such that $G_{1} \cap G_{2}$ is a complete graph. The new graph $G=G_{1} \bigoplus G_{2}$ with the vertex set $V(G)=V(G_{1}) \cup V(G_{2})$ and edge set $E(G)=E(G_{1}) \cup E(G_{2})$ is called the {\em clique sum} of $G_{1}$ and $G_{2}$ in $G_{1} \cap G_{2}$. If the cardinality of $V(G_{1}) \cap V(G_{2})$ is $k+1$, then this operation is called a $k$-clique sum of the graphs $G_{1}$ and $G_{2}$. We write $G=G_{1} \bigoplus_{\widehat{v}} G_{2}$ to indicate that $G$ is the clique sum of $G_{1}$ and $G_{2}$ and that $V(G_{1}) \cap V(G_{2})=\widehat{v}$.

\begin{thm1} \label{Exactlyone} Let $G$ be a connected graph on the vertex set $[n]$. Suppose that $G$ has exactly one cycle $C$. If $n \geq 4$ and $C$ is odd of length $3$, then ${\rm ara}(J_{G})=n-1$.

\end{thm1}

\noindent \textbf{Proof.} The graph $G$ can be written as the $0$-clique sum of the cycle $C$ and some trees. More precisely, $$G=C \bigoplus_{v_{1}} T_{1} \bigoplus_{v_{2}} \cdots \bigoplus_{v_{s}}T_{s}$$ for some vertices $v_{1},\ldots,v_{s}$ of $C$. The vertex connectivity of $G$ is 1. By Theorem \ref{Vertexconnectivity}, the inequality $n-1 \leq {\rm ara}(J_{G})$ holds. Since $G$ has exactly one cycle, we have that ${\rm card}(E(G))=n$. From Proposition \ref{Basicara}, ${\rm ara}(J_{G}) \leq n-1$, and therefore ${\rm ara}(J_{G})=n-1$. \hfill $\square$\\

Let ${\rm ht}(J_{G})$ be the height of $J_{G}$, then, we have, from the generalized Krull's principal ideal theorem, that ${\rm ht}(J_{G}) \leq {\rm ara}(J_{G})$. We say that $J_{G}$ is a {\em set-theoretic complete intersection} if ${\rm ara}(J_{G})={\rm ht}(J_{G})$.

\begin{cor1} \label{stci} Let $G$ be a connected graph on the vertex set $[n]$ with $n \geq 4$. Suppose that $G$ has exactly one cycle $C$ and its length is $3$. Then the following properties are equivalent: \begin{enumerate} \item[(a)] $J_{G}$ is unmixed, \item[(b)] $J_{G}$ is Cohen-Macaulay, \item[(c)] $J_{G}$ is a set-theoretic complete intersection, \item [(d)] $G=C \bigoplus_{v_{1}} T_{1} \bigoplus_{v_{2}} \cdots  \bigoplus_{v_{s}} T_{s}$, where $\{v_{1},\ldots,v_{s}\} \subset V(C)$, $s \geq 1$, $v_{h}$ are pairwise distinct and $T_{h}$ are paths.
\end{enumerate}
In particular, if one of the above conditions is true, then ${\rm ara}(J_{G})={\rm ht}(J_{G})=n-1$.

\end{cor1}
\noindent \textbf{Proof.} The implication (b)$\Rightarrow$(a) is well known. If $J_{G}$ is a set-theoretic complete intersection, then, from Theorem \ref{Exactlyone}, ${\rm ht}(J_{G})=n-1$ and ${\rm dim}(R/J_{G})=n+1$. Also ${\rm depth}(R/J_{G})=n+1$ by \cite[Theorem 1.1]{EHH}, so $J_{G}$ is Cohen-Macaulay, whence (c)$\Rightarrow$(b). Recall that $\mathcal{M}(G)=\{S \subset [n]: P_{S}(G) \ \textrm{is a minimal prime of} \ J_{G}\}$. If $J_{G}$ is unmixed, then every vertex $v$ of $T_{h}$, $v \neq v_{h}$, has degree at most 2. In fact, $\{v\} \in \mathcal{M}(G)$ and, if ${\rm deg}_{G}(v) \geq 3$, then by \cite[Lemma 3.1]{HTH}, one has ${\rm ht}(P_{\{v\}}(G))=n+{\rm card}(\{v\})-c(\{v\})=n+1-{\rm deg}_{G}(v) \leq n-2<n-1={\rm ht}(P_{\emptyset}(G))$, a contradiction. Moreover, $v_{h}$ has degree at most 3 for every $h$. In fact, $\{v_{h}\} \in \mathcal{M}(G)$ and, if ${\rm deg}_{G}(v_{h}) \geq 4$, then by \cite[Lemma 3.1]{HTH}, one has ${\rm ht}(P_{\{v_{h}\}}(G))=n+{\rm card}(\{v_{h}\})-c(\{v_h\})=n+1-({\rm deg}_{G}(v_{h})-1) \leq n-2<n-1={\rm ht}(P_{\emptyset}(G))$, a contradiction. Thus, (d) follows. Finally, assuming (d), $J_{G}$ is unmixed by \cite[Theorem 1.1]{EHH} and ${\rm ht}(J_{G})=n-1$. By Theorem \ref{Exactlyone}, it follows that $${\rm ara}(J_{G})=n-1={\rm ht}(J_{G}). \ \ \ \ \square$$

If $C_{1}$ and $C_{2}$ are cycles of $G$ having no common vertex, then a {\em bridge} between $C_{1}$ and $C_{2}$ is an edge $\{i,j\}$ of $G$ with $i \in V(C_{1})$ and $j \in V(C_{2})$.

\begin{prop1} \label{Bridges} Let $G$ be a connected graph on the vertex set $[n]$ with $m$ edges. Suppose that $G$ contains a subgraph $H$ consisting of two vertex disjoint odd cycles of length 3, namely $C_{1}$ and $C_{2}$, and also two bridges between the cycles $C_{1}$ and $C_{2}$. Then ${\rm ara}(J_{G}) \leq m-2$.

\end{prop1}

\noindent \textbf{Proof.} Let $E(C_{1})=\{\{1,2\},\{2,3\},\{3,1\}\}$ and $E(C_{2})=\{\{4,5\},\{5,6\},\{4,6\}\}$. Suppose first that the bridges have no common vertex. Let $e_{1}=\{1,4\}$ and $e_{2}=\{3,6\}$ be the bridges of the two cycles. Then $f_{14}^{2}$ belongs to the ideal generated by the polynomials $f_{12},f_{13},f_{14}+f_{23}$. Furthermore $f_{36}^{2}$ belongs to the ideal generated by the polynomials $f_{46},f_{56},f_{36}+f_{45}$. Thus $J_{G}$ is generated up to radical by the union of $\{f_{12},f_{13},f_{14}+f_{23},f_{46},f_{56},f_{36}+f_{45}\}$ and $\{f_{ij}|\{i,j\}\in E(G) \ \textrm{and} \ \{i,j\} \notin E(H) \}$. If the bridges have a common vertex, then without loss of generality we can assume that $e_{1}=\{1,4\}$ and $e_{2}=\{3,4\}$ are the bridges of the two cycles. Applying similar arguments as before, we deduce that ${\rm ara}(J_{G}) \leq m-2$. \hfill $\square$

\begin{ex1} {\rm Suppose that $G$ is a graph with $6$ vertices and $8$ edges consisting of two vertex disjoint odd cycles of length 3, namely $C_{1}$ and $C_{2}$, and also two vertex disjoint bridges between the cycles $C_{1}$ and $C_{2}$. Here the vertex connectivity is $l=2$. Thus ${\rm ara}(J_{G}) \geq 6$. By Proposition \ref{Bridges}, ${\rm ara}(J_{G}) \leq 6$ and therefore ${\rm ara}(J_{G})=6$.}

\end{ex1}

\begin{thm1} \label{many odd} Let $G_{k}$ be a graph containing $k$ odd cycles $C_{1},\ldots,C_{k}$ of length $3$ such that the cycles $C_{i}$ and $C_{j}$ have disjoint vertex sets, for every $1 \leq i<j \leq k$. Suppose that there exists exactly one path $P_{i,i+1}$ of length $r_{i} \geq 2$ connecting a vertex of $C_{i}$ with a vertex of $C_{i+1}$,  $1 \leq i \leq k-1$. If $G_{k}$ has no more vertices or edges, then ${\rm ara}(J_{G_k})={\rm ht}(J_{G_k})=2k+\sum_{i=1}^{k-1}r_{i}$. In particular, $J_{G_k}$ is a set-theoretic complete intersection.

\end{thm1}

\noindent \textbf{Proof.} The graph $G_{k}$ has $3k+\sum_{i=1}^{k-1}(r_{i}-1)$ vertices. Here the vertex connectivity is $l=1$, so $$2k+\sum_{i=1}^{k-1}r_{i}=3k+\sum_{i=1}^{k-1}(r_{i}-1)+1-2 \leq {\rm ara}(J_{G_k}).$$ We will prove that ${\rm ara}(J_{G_{k}}) \leq 2k+\sum_{i=1}^{k-1}r_{i}$ by induction on $k \geq 2$. Suppose that $k=2$ and let $E(C_{1})=\{\{1,2\},\{2,3\},\{1,3\}\}$, $P_{1,2}=\{\{3,4\},\{4,5\},\ldots,\{r+2,r+3\}\}$ and $C_{2}=\{\{r+3,r+4\},\{r+4,r+5\},\{r+3,r+5\}\}$. Then $J_{G_2}$ is generated up to radical by the union of $$\{f_{12}+f_{34},x_{r+2}x_{n+r+3}-x_{r+3}x_{n+r+2}+x_{r+4}x_{n+r+5}-x_{r+5}x_{n+r+4}\}$$ and $$\{f_{ij}| \{i,j\} \in E(G_{2}) \setminus \{\{1,2\},\{3,4\},\{r+2,r+3\},\{r+4,r+5\}\}\}.$$ Thus ${\rm ara}(J_{G_{2}}) \leq 4+r$. Assume that the inequality ${\rm ara}(J_{G_{k}}) \leq 2k+\sum_{i=1}^{k-1}r_{i}$ holds for $k$ and we will prove that ${\rm ara}(J_{G_{k+1}}) \leq 2(k+1)+\sum_{i=1}^{k}r_{i}$. We have that $J_{G_{k+1}}=J_{G_{k}}+J_{H}$ where $H$ is the graph consisting of the path $P_{k,k+1}$ and the cycle $C_{k+1}$. By Theorem \ref{Exactlyone}, ${\rm ara}(J_{H})=r_{k}+2$. Then, from the induction hypothesis, $${\rm ara}(J_{G_{k+1}}) \leq {\rm ara}(J_{G_{k}})+{\rm ara}(J_{H}) \leq 2k+\sum_{i=1}^{k-1}r_{i}+r_{k}+2=2(k+1)+\sum_{i=1}^{k}r_{i}.$$ Since $J_{G_k}$ is unmixed by \cite[Theorem 1.1]{EHH}, we have that $${\rm ht}(J_{G_k})={\rm card}(V(G_{k}))-1=2k+\sum_{i=1}^{k-1}r_{i}. \ \square$$

\begin{rem1} {\rm All the results presented are independent of the field $K$.\\}

\end{rem1}

\noindent \textbf{Acknowledgments}\\ 
\smallskip
\newline
The author is grateful to an anonymous referee for useful suggestions and comments that helped improve an earlier version of the manuscript. This work was supported by the Scientific and Technological Research Council of Turkey (T{\"U}BITAK) through BIDEB 2221 grant.

\end{document}